\documentclass[leqno]{article}

\usepackage{amssymb}
\usepackage{amsmath}
\usepackage{epsfig}
\newcommand{\qed}{\nobreak \ifvmode \relax \else
      \ifdim\lastskip<1.5em \hskip-\lastskip
      \hskip1.5em plus0em minus0.5em \fi \nobreak
      \vrule height0.75em width0.5em depth0.25em\fi}

\begin{document}

\title{\Large Markovian embeddings of general random strings}
\author{Manuel E. Lladser\thanks{The University of Colorado, Department of Applied Mathematics, PO Box 526 UCB, Boulder, CO 80309-0526, The United States} \thanks{e-mail: manuel.lladser@colorado.edu}}
\date{}

\maketitle

\begin{abstract}
\small\baselineskip=9pt
Let $\mathcal{A}$ be a finite set and $X$ a sequence of $\mathcal{A}$-valued random variables. We do not assume any particular correlation structure between these random variables; in particular, $X$ may be a non-Markovian sequence. An adapted embedding of $X$ is a sequence of the form $R(X_1)$, $R(X_1,X_2)$, $R(X_1,X_2,X_3)$, etc where $R$ is a transformation defined over finite length sequences. In this extended abstract we characterize a wide class of adapted embeddings of $X$ that result in a first-order homogeneous Markov chain. We show that any transformation $R$ has a unique coarsest refinement $R'$ in this class such that $R'(X_1)$, $R'(X_1,X_2)$, $R'(X_1,X_2,X_3)$, etc is Markovian. (By refinement we mean that $R'(u)=R'(v)$ implies $R(u)=R(v)$, and by coarsest refinement we mean that $R'$ is a deterministic function of any other refinement of $R$ in our class of transformations.) We propose a specific embedding that we denote as $R^X$ which is particularly amenable for analyzing the occurrence of patterns described by regular expressions in $X$. A toy example of a non-Markovian sequence of $0$'s and $1$'s is analyzed thoroughly: discrete asymptotic distributions are established for the number of occurrences of a certain regular pattern in $X_1,...,X_n$ as $n\to\infty$ whereas a Gaussian asymptotic distribution is shown to apply for another regular pattern.
\end{abstract}

%\pagenumbering{empty}

\vspace{2cm}

Full extended abstract available at:
\begin{center}
http://www.siam.org/proceedings/analco/2008/analco08.php
\end{center}

\end{document}